\documentclass[multi]{cambridge7A}
\usepackage[UKenglish]{babel}
\usepackage[longnamesfirst,sectionbib]{natbib}
\usepackage{chapterbib}
\usepackage{amsmath,amssymb,amsthm,amsfonts}
\usepackage{enumerate,url}

\setcitestyle{authoryear,round,semicolon}% alters natbib behaviour

\theoremstyle{definition}
\newtheorem{assumption}{Assumption}
\newtheorem*{ass2star}{Assumption $2^*$}

\allowdisplaybreaks[1]

\providecommand*\Index[1]{#1\index{#1}}
\providecommand*\undex[1]{} % abandoned tag
\providecommand*\Undex[1]{#1} % abandoned tag

\begin{document}
\newcommand{\BP}{\mbox{\boldmath $P$}}
\newcommand{\BQ}{\mbox{\boldmath $Q$}}
\newcommand{\Bpi}{\mbox{\boldmath $\pi$}}
\newcommand{\Bone}{\mbox{\bf 1}}
\newcommand{\Bc}{\mbox{\boldmath $c$}}
\newcommand{\Bv}{\mbox{\boldmath $v$}}
\newcommand{\Bx}{\mbox{\boldmath $x$}}
\newcommand{\BX}{\mbox{\boldmath $X$}}
\newcommand{\Bmu}{\mbox{\boldmath $\mu$}}
\newcommand{\BSigma}{\mbox{\boldmath $\Sigma$}}
\newcommand{\Balpha}{\mbox{\boldmath $\alpha$}}
\newcommand{\Bgamma}{\mbox{\boldmath $\gamma$}}
\newcommand{\var}{\mbox{var~}}

\alphafootnotes
\author[D. R. Grey]{D. R. Grey\footnotemark }
\chapter[Associated random walk and martingales]{The associated random walk
  and martingales in random walks with stationary increments}
\footnotetext[1]{Department of Probability and Statistics, Hicks Building,
  University of Sheffield, Sheffield S3 7RH; D.Grey@sheffield.ac.uk}
\arabicfootnotes
\contributor{David R. Grey \affiliation{University of Sheffield}}
\renewcommand\thesection{\arabic{section}}
\numberwithin{equation}{section}
\renewcommand\theequation{\thesection.\arabic{equation}}

\begin{abstract}
We extend the notion of the associated random walk and the Wald
martingale in random walks where the increments are independent and
identically distributed to the more general case of stationary
ergodic increments.  Examples are given where the increments are
Markovian or Gaussian, and an application in queueing is considered.
\end{abstract}

\subparagraph{AMS subject classification (MSC2010)}60G50, 60G42, 60G10, 60K25

\section{Introduction and definition}
Let $X_1$, $X_2$, \ldots\ be independent identically distributed
(i.i.d.)\ random variables with distribution function (d.f.)\ $F$ and
positive mean.  If $S_n=X_1+X_2+\cdots +X_n$ for each
$n=0$, 1, 2, \ldots\ then the process $\{S_n\} $ is called the
\index{random walk (RW)}random walk with increments $X_1$, $X_2$, \ldots. Because the
increments have
positive mean, by the strong law of large
numbers\index{strong law of large numbers (SLLN)} the random walk
will in the long run \Index{drift} upwards to infinity.

There may exist $\theta \neq 0$ such that
\[ \hat{F}(\theta ):=E(e^{-\theta X_1})=\int _{-\infty }^{\infty }
e^{-\theta x}\,dF(x)=1, \]
in which case $\theta $ is unique, and
necessarily positive because of the upward drift.  In this case,
if we define a new increment distribution by
\[ dF^*(x):=e^{-\theta x}\,dF(x), \]
then we obtain the \emph{associated random
walk}\index{random walk (RW)!associated random walk}, which has downward
drift. Because of the definition, probabilities in one random walk
may easily be expressed in terms of those in the other.  For
instance, renewal theory\index{renewal theory} in the associated random walk yields the
Cram\' er\index{Cramer, C. H.@Cram\'er, C. H.} estimate of the probability of \Index{ruin} in the original random
walk, the parameter $\theta $ determining the rate of
\Index{exponential decay} (Feller (1971), XI.7, XII.6\index{Feller, W.}; see
also Asmussen (2000))\index{Asmussen, S.}. Note
also that since $\hat{F}^*(-\theta )=\int _{-\infty }^{\infty
}e^{\theta x}\,dF^*(x)=1$ the association is a \Index{duality} relationship in
the sense that if we perform the analogous transformation on the
associated random walk with $\theta $ replaced by $-\theta $ then we
obtain the original one.

For the same value of $\theta $, it may easily be shown that
$V_n:= e^{-\theta S_n}$ defines a \Index{martingale}, known as the {\it
Wald martingale}\index{Wald, A.!Wald martingale}. This is also useful in the investigation of
hitting probabilities\index{hitting probability}.

The question arises to what extent the concepts of the associated
random walk and the Wald martingale may be generalized to the case
where the increments are no longer necessarily i.i.d.~but merely
stationary\index{stationary increments|(} and ergodic\index{ergodicity|(}.
In such cases, because of the ergodic
theorem\index{ergodicity!ergodic theorem}, the random walk still drifts upwards to infinity and so we
might still be interested in, for example, the probability of ruin
(hitting a low level).  The following is a suggested way forward. It
should be noted that, once one goes beyond the ergodic theorem, the
general stationary ergodic process behaves quite differently from
the independent-increments case.  For example, the convergence
rate\index{rate of convergence}
in the ergodic theorem may be arbitrarily slow, however many moments
may be finite.  For background, see for example Eberlein \& Taqqu
(1986).

The work below is motivated by that of Lu (1991)\index{LuZ@Lu, Z.} on
branching processes\index{branching process} in random
environments\index{random environment} (Smith \& Wilkinson
(1969)\index{Smith WL@Smith, W. L.}\index{Wilkinson, W. E.}, Athreya
\& Karlin (1971))\index{Athreya, K. B.}\index{Karlin, S.} and by a convergence result given a
straightforward proof by the author (Grey (2001)).  We obtain
generalizations of the associated random walk and the Wald
martingale, under certain assumptions.  Three applications are
considered in Section 2, to the Markov and Gaussian cases, and to a
queueing problem. Our work is also relevant to random walks in
random environments; see for example R\' ev\' esz (1990), Part
III\index{Revesz, P.@R\'ev\'esz, P.}.

To proceed, we need to make the following assumptions. In Section 2
we shall show that these assumptions are satisfied in some important
cases of interest.

\begin{assumption}There exists $\theta >0$ such that
\[ q:= \lim _{n\rightarrow \infty }E(e^{-\theta S_n}) \]
exists and is positive and finite.
\end{assumption}

This assumption is trivially satisfied in the i.i.d.~case, with
$\theta $ as identified earlier and $q=1$.  It is important to note
that since $S_n\rightarrow \infty $ with probability one, there can
be at most one value of $\theta $ satisfying the assumption.  This
is because, if such $\theta $ exists, for any positive constant $K$
\[ E(e^{-\theta S_n};S_n\leq -K)\rightarrow q \]
as the remaining contribution to the expectation tends to zero, by the
bounded convergence theorem.  From here it is easy to see that if
$0<\phi <\theta $,
\[ \limsup _{n\rightarrow \infty }E(e^{-\phi S_n};S_n\leq -K)\leq
e^{-(\theta -\phi )K}q \]
whence by similar reasoning
\[ \limsup _{n\rightarrow \infty }E(e^{-\phi S_n})\leq
e^{-(\theta -\phi )K}q \]
and so, since $K$ is arbitrary, $E(e^{-\phi S_n})\rightarrow 0$.
Similarly if $\phi >\theta $ then $E(e^{-\phi S_n})\rightarrow \infty $.

For our next assumption we extend our sequence of increments to a
doubly infinite one \ldots, $X_{-1}$, $X_0$, $X_1$, \ldots, as is always
possible with a stationary sequence (Breiman (1968), Proposition
6.5)\index{Breiman, L.}. We also define the more general partial sum
\[ S_{m,n}:=\sum _{r=m}^nX_r. \]
Let ${\cal F}_{m,n}$ denote the $\sigma $-field generated by
$\{ X_r;r=m,\ldots ,n\} $.

\begin{assumption}For all $k=1$, 2, \ldots\ and for all $B\in {\cal F}_{-k,k}$,
\[ q(B):=\lim _{m,n\rightarrow \infty }E(e^{-\theta S_{-m,n}};B) \]
exists, where $\theta $ is as defined in Assumption 1.
\end{assumption}

Again we refer immediately to the i.i.d.~case, where the assumption
is easily seen to be satisfied, the limiting operation being
essentially trivial, and we may write explicitly
\[ q(B)=E(e^{-\theta S_{-k,k}};B). \]

If we now fix $k$ and define for $m$, $n\geq k$
\[ P^*_{m,n}(B):= \frac{E(e^{-\theta S_{-m,n}};B)}{E(e^{-\theta S_{-m,n}})} \]
for $B\in {\cal F}_{-k,k}$, then clearly $P^*_{m,n}$ is a probability measure
on ${\cal F}_{-k,k}$.  Moreover if Assumptions 1 and 2 hold then
\[ P^*_{m,n}(B)\rightarrow \frac{q(B)}{q}\mbox{~~~~as~~~~}m,n\rightarrow
\infty \] for all $B\in {\cal F}_{-k,k}$.  We make use of the
result, given a straightforward proof in Grey (2001), rather simpler than for
the more general case of signed measures considered by Halmos (1950,
p.~170)\index{Halmos, P. R.}, and a special case of the
\index{Vitali, G.!Vitali--Hahn--Saks theorem}Vitali--Hahn--Saks
theorem (Dunford \& Schwartz (1958), III.7.2--4), that if $\{ P_n\} $
is a sequence of probability measures on a space $(\Omega , {\cal F})$ such
that the limit $P(A):=\lim _{n\rightarrow \infty }P_n(A)$ exists for all $A\in
{\cal F}$, then $P$ is a probability measure on $(\Omega , {\cal F})$. It
follows that
\[ P^*(B):=\frac{q(B)}{q} \]
defines a probability measure on ${\cal F}_{-k,k}$.  Since this definition
obviously does not depend upon $k$, we have consistency between different
values of $k$ and therefore a probability measure defined on
$\bigcup _{k=1}^{\infty }{\cal F}_{-k,k}$.  The
Carath\' eodory\index{Carath\'eodory, C.}
extension theorem (Durrett (1996), Appendix A.2) now ensures that $P^*$
can be extended to a probability measure defined on the whole
$\sigma $-field ${\cal F}:= {\cal F}_{-\infty ,\infty }$.

It is the probability measure $P^*$ which we use to define the
distribution of the increments of the \emph{associated random
walk}\index{random walk (RW)!associated random walk|(}.
Note that, because the original process is stationary and because of
the double-ended limiting process involved in the definition of
$P^*$, the associated process of increments is also stationary; it
would not have been possible to achieve this with a single-ended
sequence. Whether the associated process is necessarily ergodic and
whether \Index{duality} occurs is left open here;  some further remarks are
given in Section 3. Note that \Undex{ergodicity} and duality do occur in the
two special cases studied in detail in Section 2. Note also that in
the i.i.d.~case the associated random walk as defined here coincides
with the one which we have already met, since, for example, in the
discrete case
\begin{align*}
P^*(X_1=x_1,\ldots ,X_k=x_k)
  &= e^{-\theta \sum_{i=1}^kx_i}P(X_1=x_1,\ldots ,X_k=x_k)\\
  &= \prod _{i=1}^k (e^{-\theta x_i}P(X_i=x_i))
\end{align*}
as expected.

To construct our \Index{martingale}, we need to replace Assumption 2 by the following
one-sided equivalent, which again is to be read in conjunction with
Assumption 1.  Write ${\cal F}_k:={\cal F}_{1,k}$.

\begin{ass2star}For all $k=1$, 2, \ldots\ and for all $B\in {\cal F}_k$,
\[ r(B):=\lim _{n\rightarrow \infty }E(e^{-\theta S_n};B) \]
exists, where $\theta $ is as defined in Assumption 1.
\end{ass2star}

If this assumption holds, then by the aforementioned convergence
theorem, for each $k$, $r$ is a measure on ${\cal F}_k$ with total
mass $q$. Also it is absolutely continuous with respect to $P$,
since
\[ P(B)=0~~\Longrightarrow ~~\int _Be^{-\theta S_n} \,dP=0~~
\Longrightarrow ~~r(B)=0. \]
So $r$ restricted to ${\cal F}_k$ has a
Radon--Nikod\' ym\index{Radon, J.!Radon--Nikod\'ym derivative}
derivative $V_k$ with respect to $P$:
\[ \int _BV_k \,dP = r(B)\mbox{~~~~for~all~~~~}B\in {\cal F}_k, \]
where $V_k$ is ${\cal F}_k$-measurable.

But if $B\in {\cal F}_k$, then $B\in {\cal F}_{k+1}$, and
therefore
\[ \int _BV_k \,dP = r(B) = \int _BV_{k+1} \,dP\mbox{~~~~for~all~~~~}
B\in {\cal F}_k. \] So, by definition of conditional expectation,
$V_k = E(V_{k+1}|{\cal F}_k)$ almost surely, which shows that $\{
V_k\} $ is a martingale with respect to $\{ {\cal F}_k\} $.

\subparagraph{Note}In many cases it will be true that $V_k=\lim _{n\rightarrow
\infty }E(e^{-\theta S_n}|{\cal F}_k)$ almost surely, but it seems
difficult to try to use this equation as a definition of $V_k$ in general.

In the case of i.i.d.~increments it is easy to see that
$V_k=e^{-\theta S_k}$ almost surely, and so our definition
generalizes that of the Wald martingale\index{Wald, A.!Wald martingale}.

\section{Three examples}
In this section we demonstrate the existence of the associated
random walk in two important cases of interest:  stationary Markov
chain increments and stationary Gaussian increments.  In both cases,
as indeed in the simpler i.i.d.~case, certain regularity conditions
will be required.  The corresponding martingale is mentioned more
briefly in each case.  We also consider an application in queueing
theory.

\subsection{Stationary Markov chain increments}
Here we suppose that the increments $\{ X_n\} $ perform a stationary
irreducible (and therefore ergodic) \Undex{aperiodic} Markov
chain\index{Markov, A. A.!Markov chain} with countable state
space $S$.  We shall use labels such as $i$ and $j$ to represent the
actual sizes of the increments, so that they need not be integer-valued
or non-negative;  however this will not prevent us from also using them
to denote positions in vectors and matrices, since this non-standard
notation will not lead to confusion.  The associated Markov chain
constructed here has been considered in a rather more general context
by, for example, Arjas \& Speed (1973)\index{Arjas, E.}\index{Speed, T. P.}.

Let the \Index{transition matrix} of the Markov chain be $\BP=(p_{ij})$ and let
its equilibrium distribution\index{equilibrium distribution|(} be given by the column vector $\Bpi =(\pi _i)$.
Note that if $\Bone $ denotes the vector consisting entirely of ones,
then $\BP $ has Perron--Frobenius\index{Perron, O.!Perron--Frobenius eigenvalue} eigenvalue 1 with $\Bpi ^T$ and $\Bone $ as
corresponding left and right eigenvectors respectively;  also
$\Bpi ^T\Bone =1$ and $\BP ^n\rightarrow \Bone \Bpi ^T$ as $n\rightarrow
\infty $.

The regularity conditions which we impose are as follows.  For some
$\theta >0$ the matrix $\BQ $ with elements $(p_{ij}e^{-\theta j})$
has Perron--Frobenius eigenvalue 1 and corresponding left and right
eigenvectors $\Bv ^T$ (with components $(v_i)$) and $\Bc $ (with
components $(c_i)$) respectively;  also $\Bv ^T\Bc =1$ and $\BQ
^n\rightarrow \Bc \Bv ^T$ as $n\rightarrow \infty $. This
requirement is not especially stringent when the state space $S$ is
finite;  the Perron--Frobenius eigenvalue $\lambda (\theta )$ of $\BQ
$ for general $\theta $ behaves rather like the Laplace
transform\index{Laplace, P.-S.!Laplace transform}
$\hat{F}(\theta )$ in the i.i.d.~case (Lu (1991))\index{LuZ@Lu, Z.}.

We firstly show that Assumption 1 holds, with $\theta $ as defined above.
\begin{eqnarray*}
E(e^{-\theta S_n}) & = & \sum _{i_0\in S}\sum _{i_1\in S}\ldots
\sum_{i_n\in S} \pi _{i_0}p_{i_0i_1}\ldots p_{i_{n-1}i_n} e^{-\theta
(i_1+\cdots +i_n)}\\
& = & \Bpi ^T\BQ ^n\Bone ~~\rightarrow ~~\Bpi ^T\Bc \Bv
^T\Bone\mbox{~~~~as~~~~}n\rightarrow \infty .
\end{eqnarray*}

To check Assumption 2, we shall evaluate
\[ E(e^{-\theta S_{-m,n}};X_{-k}=i_{-k},\ldots ,X_k=i_k) \]
for given $k<m$, $n$ and $i_{-k}$, \ldots, $i_k\in S$.  By the Markov
property\index{Markov, A. A.!Markov property}, this may be written as the product of the three factors $
e^{-\theta (i_{-k}+\cdots +i_k)}\allowbreak P(X_{-k}=i_{-k},\ldots ,X_k=i_k)$
together with $E(e^{-\theta S_{k+1,n}}|X_k=i_k)$ and\break $E(e^{-\theta
S_{-m,-k-1}}|X_{-k}=i_{-k})$. The first factor may be written
\[ e^{-\theta (i_{-k}+\cdots +i_k)}\pi _{i_{-k}}p_{i_{-k}i_{-k+1}}\ldots
p_{i_{k-1}i_k}. \]
The second factor may be written
\[ \sum _{i_{k+1}\in S}\ldots \sum _{i_n\in S}e^{-\theta (i_{k+1}+\cdots
+i_n)}p_{i_ki_{k+1}}\ldots p_{i_{n-1}i_n}. \]
Because the reverse Markov chain\index{reverse Markov chain} has transition probabilities
$\pi _jp_{ji}/\pi _i$, the third factor may be written
\[ \sum _{i_{-k-1}\in S}\!\cdots\! \sum _{i_{-m}\in S}e^{-\theta (i_{-k-1}+\cdots
+i_{-m})}\frac{\pi _{i_{-k-1}}}{\pi _{i_{-k}}}p_{i_{-k-1}i_{-k}}\cdots
\frac{\pi _{i_{-m}}}{\pi _{i_{-m+1}}}p_{i_{-m}i_{-m+1}}. \]
The second factor is seen to be the $i_k$ component of the vector
$\BQ ^{n-k}\Bone $ and so converges to $c_{i_k}\Bv ^T\Bone $ as
$n\rightarrow \infty $.

Writing $\mu _i:=\pi _ie^{-\theta i}$ for each $i\in S$ and letting
$\Bmu $ be the corresponding vector, after cancellation and
rearrangement the third factor is seen to be $\mu _{i_{-k}}^{-1}$
times the $i_{-k}$ component of the vector $\Bmu ^T\BQ ^{m-k}$ and
so converges to $\mu _{i_{-k}}^{-1} v_{i_{-k}}\Bmu ^T\Bc $ as
$m\rightarrow \infty $.  Note that
\begin{align*}
\Bmu ^T\Bc&=\sum _{j\in S}\pi_je^{-\theta j}c_j\\
  &=\sum_{j\in S}\sum_{i\in S}\pi_ip_{ij}e^{-\theta j}c_j
  =\sum_{i\in S}\pi_i\sum_{j\in S}p_{ij}e^{-\theta j}c_j
  =\sum_{i\in S}\pi_ic_i
  =\Bpi^T\Bc.
\end{align*}
Using this fact and putting all the preceding results together,
after cancellation we see that $P^*(X_{-k}=i_{-k},\ldots ,X_k=i_k)$
exists and is equal to
\[ e^{-\theta (i_{-k+1}+\cdots +i_k)}p_{i_{-k}i_{-k+1}}\ldots
p_{i_{k-1}i_k}c_{i_k}v_{i_{-k}}. \]
Writing $p^*_{ij}:=p_{ij}e^{-\theta j}c_j/c_i$ and $\pi ^*_i:=c_iv_i$ for
each $i$, $j\in S$, the above may also be written
\[ P^*(X_{-k}=i_{-k},\ldots ,X_k=i_k) = \pi ^*_{i_{-k}}p^*_{i_{-k}i_{-k+1}}
\ldots p^*_{i_{k-1}i_k}. \]
It is a routine matter to check that the numbers $(p^*_{ij})$ form the
transition probabilities of a Markov chain and that $(\pi ^*_i)$ is an
equilibrium distribution\index{equilibrium distribution|)} for it.  We have thus established that the
associated random walk exists and its increments perform a stationary
Markov chain, which is also obviously irreducible and \Undex{aperiodic} like
the original.  Duality\index{duality} is left as an exercise.

The martingale\index{martingale|(} may be constructed similarly.  Letting $B=\{ X_1=i_1,\ldots ,
X_k=i_k\} $ it may be calculated that
\[ E(e^{-\theta S_n};B)\rightarrow P(B)e^{-\theta (i_1+\cdots +i_k)}
c_{i_k}\Bv ^T\Bone\mbox{~~~~as~~~~}n\rightarrow \infty \] and so
Assumption $2^*$ holds;  moreover, since $B$ is an atom of the
$\sigma $-field ${\cal F}_k$ it follows that $V_k = e^{-\theta
S_k}c_{X_k}\Bv ^T\Bone $. We may take $\Bv ^T\Bone =1$ since $\Bv $
and $\Bc $ have so far only been scaled relative to each other. This
gives the martingale\index{martingale|)} $V_k = c_{X_k}e^{-\theta S_k}$ which has also
been used by Lu (1991)\index{LuZ@Lu, Z.}.

\subsection{Stationary Gaussian increments}
Now let $\{ X_n\} $ be a stationary Gaussian
process\index{Gauss, J. C. F.!Gaussian process} in which each $X_n$
has normal distribution with mean $\mu >0$ and variance $\sigma ^2 >0$.
For each $r=1$, 2, \ldots\ let $\rho _r$ be the \Index{correlation} coefficient
between $X_n$ and $X_{n+r}$ for any $n$.  These parameters completely
determine the behaviour of the process, since the joint distribution
of any finite collection of the $X_n$ is multivariate normal.

The regularity condition we need here is that
\[ \sum _{r=1}^{\infty }r|\rho _r| <\infty . \]
This is an \Index{asymptotic independence} property more than sufficient for
\Undex{ergodicity}, and one which is easily satisfied by commonly studied
processes such as autoregressive\index{autoregressive process} and moving
average processes\index{moving average process}.

Under this condition, let $R:=\sum _{r=1}^{\infty }\rho _r$ and let
$S:=\sum _{r=1}^{\infty }r\rho _r$; these will both be finite.
Below we shall see that $R\geq -\frac{1}{2}$ necessarily, and that we
need to exclude the extreme case $R=-\frac{1}{2}$.

We firstly find $\theta $ such that Assumption 1 is satisfied.
Since $S_n=X_1+\cdots +X_n$ has a normal distribution with mean $n\mu $
and variance $\sigma ^2[n+2\sum _{r=1}^{n-1}(n-r)\rho _r]$, by the
standard formula for the Laplace transform\index{Laplace, P.-S.!Laplace transform} of the normal
distribution we have that
\[ E(e^{-\theta S_n})=\exp \left\{ -n\mu \theta +
{\textstyle \frac12 }\sigma ^2 [n+2\sum _{r=1}^{n-1}(n-r)\rho
_r]\theta ^2 \right\} . \] Under our regularity condition
\[ \sum _{r=1}^{n-1}(n-r)\rho _r = nR-S+o(1)\mbox{~~~~as~~~~}
n\rightarrow \infty \]
and so in particular
\[ \var S_n = \sigma ^2(n[1+2R]-2S)+o(1)\mbox{~~~~as~~~~}
n\rightarrow \infty, \] whence $1+2R\geq 0$.  It is possible to
construct examples with $1+2R=0$ (such as $X_n:=\mu +Z_n-Z_{n-1}$
where $\{ Z_n\} $ are i.i.d. $N(0,\frac{1}{2}\sigma ^2)$ random
variables) but if we exclude this rather extreme case then we see
that convergence of $E(e^{-\theta S_n})$ to a positive limit will
occur if and only if
\[ -\mu \theta + {\textstyle \frac12 }\sigma ^2[1+2R]\theta ^2 = 0. \]
This yields
\[ \theta = \frac{2\mu }{\sigma ^2[1+2R]} \]
and then it is easy to compute that for this value of $\theta $
\[ E(e^{-\theta S_n})\rightarrow \exp \left\{ -\frac{4\mu ^2S}
{\sigma ^2[1+2R]^2} \right\}\mbox{~~~~as~~~~}n\rightarrow \infty . \]

We turn to Assumption 2.  Fix $k$, and take $\theta $ as just
identified. For $m$, $n>k$ it is evidently relevant to look at the
distribution of $Y:=S_{-m,-k-1}+S_{k+1,n}$ conditional on
$X_{-k}=x_{-k}$, \ldots, $X_k=x_k$, or $\BX =\Bx $ say.  If the
unconditional distribution of $Y$ is denoted by $N(\nu ,\tau ^2)$,
the vector of covariances of $Y$ with the components of $\BX $ is
denoted by $\Bv $, and the mean and covariance matrix of $\BX $ are
denoted by $\Bmu $ and $\BSigma $ respectively, then by multivariate
normal theory (Mardia, Kent \& Bibby (1979), Theorem
3.2.4)\index{Mardia, K. V.}\index{Kent, J. T.}\index{Bibby, J. M.}, the
conditional distribution is $N(\nu +\Bv ^T\BSigma ^{-1}(\Bx -\Bmu
),\tau ^2-\Bv ^T\BSigma ^{-1}\Bv ).$ (For ease of notation, we
suppress the dependence of $Y$, $\nu$, $\tau ^2$ and $\Bv $ on $m$ and
$n$.)

A typical component of $\Bv $ is of the form
\[ \sigma ^2\left(\sum _{r=i+k+1}^{i+m} \rho _r+\sum _{r=k+1-i}^{n-i} \rho _r\right) \]
for some $i$, and so $\Bv $ converges to a finite limit as
$m$, $n\rightarrow \infty $.  Also $\nu =(m+n-2k)\mu $.
Then, since
\[ \var S_{-m,n} = \var (\Bone ^T\BX +Y) = \Bone ^T\BSigma \Bone +
2\Bone ^T\Bv +\tau ^2, \] we can use the estimate of $\var S_n$
obtained in checking Assumption 1 to deduce that
\[ \tau ^2 -\sigma ^2[1+2R](m+n+1) \]
converges to a finite limit as $m$, $n\rightarrow \infty $. Putting these results
together we see that
\[ E(e^{-\theta Y}|\BX =\Bx ) = \exp \left\{ -[\nu +\Bv ^T\BSigma ^{-1}
(\Bx -\Bmu )]\theta + {\textstyle \frac12 }[\tau ^2-\Bv ^T\BSigma
^{-1}\Bv]\theta ^2 \right\} \] converges to a positive limit as
$m$, $n\rightarrow \infty $, since because of the value of $\theta $
the difference between the large terms in $\nu $ and $\tau ^2$
converges to a finite limit, and the other terms also converge to
finite limits. We may denote the limit in the above by $\exp
(\Balpha ^T\Bx +\beta )$ for some constants $\Balpha $ and $\beta $.
It is then easy to see that for any $B\in {\cal F}_{-k,k}$,
\[ E(e^{-\theta S_{-m,n}};B) \rightarrow E(\exp \{-\theta \Bone ^T\BX +\Balpha ^T
\BX +\beta \};B)\quad\mbox{as}\quad m,\,n\rightarrow\infty,\] and we
have established Assumption 2.  It is not hard in this case to see
that the associated random walk has the same covariance structure as
the original, but downward drift $-\mu $.

Calculations similar to the above may be used to check that Assumption
$2^*$ holds in this case also, and that the associated \Index{martingale} is of the
form $V_k = \exp (\Bgamma ^T\BX +\delta )$ for some $\Bgamma$, $\delta $,
where now $\BX=(X_1,\ldots ,X_k)^T$.\index{stationary increments|)}\index{random walk (RW)!associated random walk|)}

\subsection{A queueing application}
Suppose we have a $G/GI/1$ queue\index{queue|(} in which the inter-arrival times
$\{ T_n\} $ form a stationary ergodic\undex{ergodicity} sequence and the independent
service times $\{ U_n\} $ are i.i.d.~with $ET_n > EU_n > 0$.  Then
the waiting times $\{ W_n\} $ satisfy
\[ W_{n+1} = (W_n + U_n - T_n)^+ \]
and it follows by a standard argument, dating back to Lindley
(1952)\index{Lindley, D. V.}
in the case of i.i.d.~inter-arrival times and exploited, among
others, by Kingman (1964)\index{Kingman, J. F. C.}, that $W_n$ has an
\Index{equilibrium distribution} which is the same as the distribution of minus the
all-time minimum of an unrestricted
\index{random walk (RW)}random walk started at zero in
state zero, with increments $X_n:=T_{-n}-U_{-n}$.  The tail of this
distribution is therefore intimately related to the probability of
\Index{ruin} in this random walk, and, in particular, the parameter $\theta
$, if it exists, has an important part to play.

The simplest example is the $M/M/1$ queue where the $T_n$ are
independent exponential with parameter $\lambda $ and the $U_n$ are
exponential with parameter $\mu $, where $\mu > \lambda > 0$. In
this case the $X_n$ are i.i.d. and $E(e^{-\theta X_n}) = \lambda \mu
/ \{ (\lambda + \theta )(\mu - \theta ) \} $, which is easily seen
to be equal to one when $\theta = \mu - \lambda $.  So the key
parameter $\theta $ depends upon both the arrival rate and the
service rate in a simple and obvious way.

As another example, suppose that the $U_n$ have some arbitrary
distribution with Laplace transform\index{Laplace, P.-S.!Laplace transform} $\phi (\theta )=E(e^{-\theta
U_n})$, and that there is a regular \Index{appointments system} such that
customer $n$ arrives at clock time $\lambda ^{-1}n + \epsilon _n$,
where $\{ \epsilon _n\} $ is a sequence of i.i.d.~errors with
Laplace transform $\psi (\theta ) = E(e^{-\theta \epsilon _n})$.  In
this case it is possible to compute
\[ E(e^{-\theta S_n}) = [\phi (-\theta )]^n \exp (-\lambda ^{-1}n\theta
) \psi (\theta ) \psi (-\theta ) \] and so the key parameter $\theta
$ satisfies the equation
\[ \phi (-\theta ) \exp (-\lambda ^{-1}\theta ) = 1 \]
which does not involve the distribution of the $\epsilon _n$.  By
considering the special case $\phi (\theta )=\mu /(\mu +\theta )$,
so that the $U_n$ are exponentially distributed, and the mean
inter-arrival time and mean service time are $\lambda ^{-1}$ and
$\mu ^{-1}$ respectively as in the previous $M/M/1$ example, it is
possible to compare an appointments system with random
(Poisson)\index{Poisson, S. D.!Poisson process}
arrivals.  For the value $\theta =\mu -\lambda $ found in the case
of random arrivals, we may compute $\phi (-\theta ) \exp (-\lambda
^{-1}\theta )=(\mu /\lambda )\exp (1-(\mu /\lambda ))<1$, and so the
actual value of the key parameter $\theta $ for the appointments
system must be larger than for random arrivals.  This suggests a
thinner tail for the equilibrium waiting time distribution, and a
more efficient system.\index{queue|)}

\section{Some remarks on duality and asymptotic independence}
If we wish for duality\index{duality|(} to occur, then, replacing $\theta $ by $-\theta $
and denoting expectation with respect to $P^*$ by $E^*$, we require for
$B\in {\cal F}_{-k,k}$ that
\[ \frac{E^*(e^{\theta S_{-m,n}};B)}{E^*(e^{\theta S_{-m,n}})}\rightarrow
P(B)\mbox{~~~~as~~~~}m,n\rightarrow \infty . \] Now, denoting
$E(e^{-\theta S_{-r,s}})$ by $q_{r,s}$ and noting that $P^*_{r,s}$
has Radon--Nikod\' ym derivative\index{Radon, J.!Radon--Nikod\'ym derivative} $q_{r,s}^{-1}e^{-\theta S_{-r,s}}$
with respect to $P$, we have that
\begin{align*}
E^*(e^{\theta S_{-m,n}};B)&=\int_B e^{\theta S_{-m,n}}\,dP^*\\
&=\lim_{r,s\rightarrow\infty}\int_B e^{\theta S_{-m,n}}\,dP^*_{r,s}\\
&=\lim_{r,s\rightarrow\infty}
    \int_B e^{\theta S_{-m,n}}q_{r,s}^{-1}e^{-\theta S_{-r,s}}\,dP\\
&=q^{-1}\lim_{r,s\rightarrow\infty}
    \int_B e^{-\theta S_{-r,-m-1}}e^{-\theta S_{n+1,s}}\,dP.
\end{align*}
For the required convergence to occur,
it seems therefore that for large $m$ and $n$ there must be
approximate independence between $S_{-r,-m-1}$, $S_{n+1,s}$ and the
event $B$, so that we can say that the above integral is
approximately $P(B) E(e^{-\theta S_{-r,-m-1}}) E(e^{-\theta
S_{n+1,s}})$ which converges to $P(B)q^2$ as $r$, $s\rightarrow \infty
$.  Hence under these circumstances
\[ \frac{E^*(e^{\theta S_{-m,n}};B)}{E^*(e^{\theta S_{-m,n}})} \sim
\frac{q^{-1}P(B)q^2}{q^{-1}q^2}\rightarrow P(B)\mbox{~~~~as~~~~}
m,n\rightarrow \infty . \] The \Index{asymptotic independence} is a kind
of \Index{mixing} condition which is already stronger than
\index{ergodicity|)}ergodicity,
suggesting that the latter is not the most appropriate property to
be considering in this context.  See Bradley (2005)\index{Bradley, R. C.} on mixing
conditions.\index{duality|)}

\paragraph{Acknowledgements}
The author would like to thank Nick Bingham and an anonymous referee
for helpful remarks on the presentation of this paper.

\end{document}